\documentclass[a4paper,12pt]{article}
\usepackage{amssymb,amsfonts,amsmath,amsthm}
\usepackage[cp1251]{inputenc}
\usepackage[english]{babel}
\title{Local maxima of energy of point configurations on sphere}
\author{Danylo Radchenko}
\newtheorem{thm}{Theorem}

\begin{document}
\maketitle

\smallskip

\begin{abstract}
We prove that certain energy functionals of point configurations on sphere have no local maxima.\end{abstract}

In this short note we will prove that certain energy functionals on sphere have no local maxima. This partially answers the question that Professor E.Saff asked on the conference "Optimal Configurations on the Sphere and Other Manifolds" at Vanderbilt University in 2010.

For $d\in \mathbb{N}$, denote by $S^d$ a unit sphere in $\mathbb{R}^{d+1}$. For $\alpha>0$ and any configuration of $N\geq 2$ distinct points $x_1,x_2,\ldots,x_N\in S^d$ consider the following energy functional
$$ E_{\alpha}(x_1,x_2,\ldots,x_N)=\sum_{i\neq j}\frac{1}{\|x_i-x_j\|^{\alpha}},$$
where $\|\cdot\|$ is the Euclidean norm on $\mathbb{R}^{d+1}$.

 For $d=2,\alpha=1$ this functional has a physical interpretation as the electrostatic potential energy of a system containing $N$ equally charged particles on the sphere.

 The problem of finding configurations which minimize these functionals has strong connections to the problem of finding uniformly distributed collections of points on sphere (see Saff, Kuijlaars \cite{SK}), as well as to the problem of finding optimal spherical codes (see Cohn, Kumar \cite{CK}).

 It is easy to see that for each $d\in \mathbb{N}, N\geq 2$ and $\alpha>0$ there exists a configuration of $N$ points on sphere $S^d$ at which $E_{\alpha}$ has a local (and even global) minimum. However it is not entirely clear whether $E_{\alpha}$ can have local maxima. We prove the following theorem, which states that for large $\alpha$ this never happens.

\begin{thm}
 For $\alpha\geq d-2$ the functional $E_{\alpha}(x_1,x_2,\ldots,x_N)$ has no local maxima.
\end{thm}

\emph{Proof:} For convenience, we scale the energy by a factor of $2^{\alpha/2}$.  Put $r=\alpha/2$ and denote $g_r(t)=(1-t)^{-r}$. Then we have
 $$\frac{2^r}{\|x_i-x_j\|^{\alpha}}=g_r(\langle x_i,x_j\rangle),$$
where $\langle\cdot,\cdot\rangle$ is the usual inner product on $\mathbb{R}^{d+1}$. Then we can write the energy functional as
 $$ E_{\alpha}(x_1,\ldots,x_N)=\sum_{i\neq j} g_r(\langle x_i,x_j\rangle).$$

Introduce arbitrary vectors $h_1,h_2,\ldots,h_N$ orthogonal to corresponding $x_i$, (i.e. $\langle x_i,h_i\rangle = 0$) and consider the function $f:\mathbb{R}\rightarrow\mathbb{R}$ defined by
$$ f(t)=E_{\alpha}\Big(\frac{x_1+th_1}{\|x_1+th_1\|},\ldots,\frac{x_N+th_N}{\|x_N+th_N\|}\Big).$$
If $E_{\alpha}$ attains a local maximum at $x_1,\ldots,x_N$, then we must have $f'(0)=0$ and $f''(0)\leq 0$. After some elementary calculations one can verify the following identity
\begin{multline}
f''(0)=\sum_{i\neq j} \Big[g_r''(\langle x_i,x_j\rangle)(\langle x_i,h_j\rangle+\langle x_j,h_i\rangle)^2+ \\
g_r'(\langle x_i,x_j\rangle)\big(2\langle h_i,h_j\rangle-(\|h_i\|^2+\|h_j\|^2)\langle x_i,x_j\rangle\big)\Big].
\label{eq:n1}
\end{multline}
Therefore, in order to prove that our energy has no local maxima it is sufficient to find $h_i$ such that (\ref{eq:n1}) is strictly positive. To do so, take $h_2=h_3=\ldots=h_N=0$ and $h_1=h$, where $\|h\|=1$. Then $f''(0)/2$ is equal to
\begin{equation}
 \sum_{j=2}^{N}\big[g_r''(\langle x_1,x_j\rangle)\langle x_j,h\rangle^2-g_r'(\langle x_1,x_j\rangle)\langle x_1,x_j\rangle \big].\label{eq:n2}
\end{equation}

Suppose that (\ref{eq:n2}) is $\leq 0$ for all $h$ orthogonal to $x_1$. Then the average value of (\ref{eq:n2}) over all such $h$ is also $\leq 0$. More specifically, let $H=\{h\in S^d: \langle x_1,h\rangle=0\}$, then $H$ is a $(d-1)$-dimensional sphere, and we take $\mu_{d-1}$ to be the normalized Lebesgue measure on $H$. We have
$$\int_{H}\langle x_j,h\rangle^2d\mu_{d-1}(h)=\int_{H}\langle x_j-x_1\langle x_1,x_j\rangle,h\rangle^2d\mu_{d-1}(h)=\frac{1-\langle x_j,x_1\rangle^2}{d},$$
because $x_j'=x_j-x_1\langle x_1,x_j\rangle$ belongs to $H$ and $\|x_j'\|^2=1-\langle x_j,x_1\rangle^2$. Therefore, integrating (\ref{eq:n2}) over $H$ with respect to $h$ gives us
\begin{equation}
 \sum_{j=2}^{N}\big(g_r''(\langle x_1,x_j\rangle)\frac{1-\langle x_1,x_j\rangle^2}{d}-g_r'(\langle x_1,x_j\rangle)\langle x_1,x_j\rangle\big) \leq 0. \label{eq:n3}
\end{equation}
After substituting $g_r(t)=(1-t)^{-r}$ into (\ref{eq:n3}) we get
$$
 \sum_{j=2}^{N}\Big(\frac{r(r+1)(1+\langle x_1,x_j\rangle)}{d(1-\langle x_1,x_j\rangle)^{r+1}}-\frac{r\langle x_1,x_j\rangle}{(1-\langle x_1,x_j\rangle)^{r+1}}\Big)\leq 0,
$$
or, equivalently,
\begin{equation}
 \sum_{j=2}^{N}r\frac{(r+1)+(r+1-d)\langle x_1,x_j\rangle}{d(1-\langle x_1,x_j\rangle)^{r+1}}\leq 0. \label{eq:n4}
\end{equation}
Since $\alpha \geq d-2$, we have $|r+1-d|\leq r+1$ and hence every term on the left of (\ref{eq:n4}) is nonnegative. In fact, we have $\langle x_1,x_j\rangle<1$, so every term is strictly positive, and therefore the sum on the left of (\ref{eq:n4}) must be strictly positive. This contradiction concludes the proof. \qed

\end{document}